\documentclass[11pt]{extarticle}
\usepackage[english]{babel}
\usepackage{graphicx}
\usepackage{framed}
\usepackage[top=1 in,bottom=1in, left=1 in, right=1 in]{geometry}
\usepackage[mathlines]{lineno}

\usepackage[utf8]{inputenc} 
\usepackage[T1]{fontenc}    
\usepackage{hyperref}       
\usepackage{url}            
\usepackage{booktabs}       
\usepackage{amsfonts}       
\usepackage{amsmath}
\usepackage{amsthm}
\usepackage{amssymb}
\usepackage{mathtools}
\usepackage{commath}
\usepackage{mathrsfs}
\usepackage{bbm}
\usepackage{nicefrac}       
\usepackage{microtype}      
\usepackage{lipsum}
\usepackage{graphicx}
\graphicspath{ {./images/} }
\usepackage{xcolor}
\usepackage{appendix}
\usepackage{enumitem}   

\newtheorem{theorem}{Theorem}
\newtheorem{proposition}{Proposition}
\newtheorem{lemma}{Lemma}
\newtheorem*{lemma*}{Lemma}
\newtheorem*{notations*}{General notations}

\newtheorem*{remark*}{\textit{Remark}}
\newtheorem*{facts*}{\textit{Facts}}

\newtheorem{assumption}{Assumption}

\newtheorem*{definition*}{Definition}

\newcommand{\m}{\mathbf{m}}
\newcommand{\M}{\mathbf{M}}

\newcommand{\R}{\mathbb{R}}
\newcommand{\Rd}{\mathbb{R}^d}
\newcommand{\E}{\mathbb{E}}

\title{Mean-field limit of Age and Leaky Memory dependent Hawkes processes}

\author{
 Valentin Schmutz\footnote{Brain Mind Institute, École Polytechnique Fédérale de Lausanne (EPFL), 1015 Lausanne, Switzerland (\texttt{valentin.schmutz@epfl.ch}).}
}
\date{}
\begin{document}
\maketitle

\begin{abstract}
We propose a mean-field model of interacting point processes where each process has a memory of the time elapsed since its last event (age) and its recent past (leaky memory), generalizing Age-dependent Hawkes processes. The model is motivated by interacting nonlinear Hawkes processes with Markovian self-interaction and networks of spiking neurons with adaptation and short-term synaptic plasticity. 

By proving propagation of chaos and using a path integral representation for the law of the limit process, we show that, in the mean-field limit, the empirical measure of the system follows a multidimensional nonlocal transport equation. 
\end{abstract}

    


\section{Introduction}
The dynamics of many interacting particle systems can be approximated, when the size of the system tends to infinity, by a partial differential equation (PDE) \cite{KipLan99}. This not only links microscopic and macroscopic scales but also stochastic and deterministic models. For mean-field models, one can prove this type of results by exploiting the propagation of chaos phenomenon, i.e. for \textit{i.i.d.} initial conditions, particles become asymptotically independent in the mean-field limit \cite{Szn91,Mel96}.

Propagation of chaos arguments have been applied to the study of interacting point processes \cite{DelFou16, GalLoe16, Loe17, MehSch20}. This has been particularly important for the field of theoretical neuroscience as it has provided a rigorous footing to the \textit{population density} formalism, where the dynamics of a population of neurons is described by a PDE (see \cite[Part~III]{GerKis14}). An example of population density equation is the refractory density (or age-structured) equation \cite{GerHem92,Ger95,GerKis02, PakPer10, PakPer13, GerKis14,SchChi19}, which has recently be proved to be exact in the mean-field limit \cite{DemGal15, Qui16, Che17}.

The models considered in \cite{DemGal15, Qui16, Che17} all assume that the point processes are `renewal' (in some loose sense), i.e. each process has a memory of its past that is restricted to the time elapsed since its last event. The fact that, in the `renewal' case, the mean-field limit can by characterized by relatively simple deterministic equations has long been recognized in theoretical neuroscience and has led to a large body of work \cite{WilCow72, Ger95, Ger00, SchChi19}. In contrast, the case where point processes are not `renewal' is much less understood. In particular, even though some heuristic population density equations have been proposed for the `non-renewal' case \cite{MulBue07,ToyRad09,PakPer14}, their exactness in the mean-field limit has not been proved. The aim of this work is therefore to propose a general framework for relating interacting `non-renewal' point processes with PDEs, in the mean-field limit. This framework relies on the definition of an abstract interacting point process model, which generalizes Age-dependent Hawkes processes \cite{Che17,RaaDit20}.

\subsection{Interacting Age and Leaky Memory dependent Hawkes processes}
Consider a system of $N$ interacting point-processes, interacting through a common variable $X^N_t$. Each point process $i$ is associated with $1+d$ variables (for $d$ a positive integer): an \textit{age} variable $A^{i,N}_t$ which represents the time elapsed since the last event of process $i$ and a $d$-dimensional vector of \textit{leaky memory} variables $\M^{i,N}_t$ which models the effect of the recent past of process $i$. The point process $i$ has stochastic intensity $(f(A^{i,N}_{t-},\M^{i,N}_{t-}, X^N_{t-}))_{t\in\R_+}$ where $f: \R_+ \times \mathbb{R}^d \times \R \to \mathbb{R}_+$ is the intensity function. Intuitively, this means that if we write $(Z^{i,N}_t)_{t\in\R_+}$ the counting process associated with the point process $i$, the instantaneous probability for $Z^{i,N}$ to jump in $]t,t+dt]$, given the past $\mathcal{F}_{t}$, is
\begin{linenomath}\begin{equation*}
    \mathbb{P}(Z^{i,N}_{t+dt}>Z^{i,N}_{t}|\mathcal{F}_{t}) = f(A^{i,N}_{t},\M^{i,N}_{t}, X^N_{t})dt.
\end{equation*}\end{linenomath}

Between events (jumps) of process $i$, the \textit{age} variable $A^{i,N}_t$ grows linearly with time whereas the \textit{leaky memory} variables $\M^{i,N}_t$ drift following the vector field $b:\Rd \to \Rd$.

At each event of process $i$, its \textit{age} $A^{i,N}_t$ is reset to $0$ and its \textit{leaky memory} $\M^{i,N}_t$ jumps to $\M^{i,N}_t + \Gamma(\M^{i,N}_t)$, where $\Gamma: \R^d \to \R^d$ is the \textit{jump function}. The fact that the variables $\M^{i,N}_t$ are not reset to a fixed value at each event allows them to accumulate the effect of successive events.

Finally, the time-dependent effect of an event of point process $i$ on point process $j$ is determined by the interaction function $h: \R_+ \times \R_+ \times \R^d \to \R$ which depends on $A^{j,N}_t$ and $\M^{j,N}_t$. Since the function $h$ is the same for all $i$ and $j$, the interaction is said to be of mean-field type. 

The model can be described by a system of stochastic integral equations: for $i = 1, \dots, N$, 
\begin{subequations}\label{eq:ALM-Model}
\begin{linenomath}\begin{align}
    A^{i,N}_t &= A^i_0 + t - \int_0^t A^{i,N}_{s-}dZ^{i,N}_s, \\
    \M^{i,N}_t &= \M^i_0 + \int_0^t b(\M^{i,N}_s)ds + \int_0^t \Gamma(\M^{i,N}_{s-})dZ^{i,N}_s, \label{eq:ALM-Model-M}\\
    Z^{i,N}_t &= \int_{[0,t]\times\R_+}\mathbbm{1}_{z\leq f(A^{i,N}_{s-}, \M^{i,N}_{s-}, X^{N}_{s-})}\pi^i(ds,dz),
\end{align}\end{linenomath}
with
\begin{linenomath}\begin{equation}\label{eq:ALM_X}
X^N_t = \frac{1}{N}\sum_{j=1}^N H^j_t + \frac{1}{N}\sum_{j=1}^N\int_0^t h(t-s,A^{j,N}_{s-}, \M^{j,N}_{s-})dZ^{j,N}_s.
\end{equation}\end{linenomath}
\end{subequations}
The collection $\{\pi^i\}_{i\in\mathbb{N}^*}$ is a sequence of independent Poisson random measures on $\R_+\times\R_+$ with Lebesgue intensity measure. We work on the filtered probability space $(\Omega, \mathcal{F}, (\mathcal{F}_t)_{t\in\R_+}, \mathbb{P})$ where $\{\pi^i\}_{i\in\mathbb{N}^*}$ is independent of $\mathcal{F}_0$ and $\mathcal{F}_t = \mathcal{F}_0 \cup \sigma\left(\{\pi^i([0,t],B)\}_{i\in\mathbb{N}^*, B\in\mathcal{B}(\R_+)}\right)$. For all $i\in \mathbb{N}^*$, $A^i_0$ and $\M^i_0$ are $\mathcal{F}_0$-measurable random variables taking values in $\R_+$ and $\Rd$ respectively and $(H^i_t)_{t\in\R_+}$ is a $\mathcal{F}_0$-measurable $\mathcal{C}(\R_+)$ random function. The $1/N$ scaling in~\eqref{eq:ALM_X} will allow us to take the mean-field limit $N\to\infty$.

If $f$ does not depend on the leaky memory variables $\M$ and $h$ does not depend on the age $A$ nor $\M$, \eqref{eq:ALM-Model} reduces to a system of interacting Age-dependent Hawkes processes \cite{Che17,RaaDit20}. If, in addition, $f$ does not depend on $A$, the model further reduces to a mean-field system of interacting nonlinear Hawkes processes (with vanishing self-interaction) \cite{DelFou16,DitLoe17}. The model~\eqref{eq:ALM-Model} has two motivations: first, it is general enough to encompass several concrete examples from the theory of nonlinear Hawkes processes and neuroscience (see below); second, its mean-field limit can be characterized by a PDE.

\subsection{Motivating examples}

Hawkes processes \cite{Haw71} provide a flexible and intuitive model for point processes with dependence on the past. They have found applications in finance \cite{BacMas15,Haw18}, seismology \cite{Oga99}, social systems \cite{CraSor08}, genomics \cite{ReySch10} and  neuroscience \cite{TruEde05,PilShl08,ReyRiv13,GerEde13,Tru16,LamTul18}, among other fields. Neuroscience research has mainly focused on \textit{nonlinear} Hawkes processes \cite{BreMas96} since they are closely related to well-established neuron models such as the Spike Response Model \cite{Ger00, JolRau06, GerKis02, GerKis14} and the Recursive Linear-Nonlinear Poisson Model \cite{PilShl08}, both variations of Generalized Linear Models (see \cite{GerKis14} Part II and
references therein). However, the models differ from the nonlinear Hawkes processes considered in \cite{DelFou16,DitLoe17} in that, even when $N$ is large, self-interaction (the effect of process $i$ on itself) does not vanish. Self-interaction vanishes in \cite{DelFou16,DitLoe17} because it is scaled by $1/N$. Let us now consider the case where self-interaction $\mathfrak{h}$ can be different from hetero-interaction $h$ (the effect of a process on the other processes) and only hetero-interaction is scaled by $1/N$: for $i = 1, \dots, N$,
\begin{subequations} \label{eq:Nd-hawkes-general}
\begin{linenomath}\begin{align}
    Z^{i,N}_t &= \int_{[0,t]\times \R_+} \mathbbm{1}_{\{z\leq f(X^{i,N}_{s-})\}}\pi^i(ds,dz), \\
    X^{i,N}_t(i) &= \mathfrak{H}^i_t + \frac{1}{N}\sum_{j\neq i} H^j_t + \int_0^t \mathfrak{h}(t-s)dZ^{i,N}_s + \frac{1}{N}\sum_{j\neq i}\int_0^t h(t-s)dZ^{j,N}_s,
\end{align}\end{linenomath}
\end{subequations}
where $f$ is a monotonically increasing function and $\{(\mathfrak{H}_t(i))_{t\in\R_+}\}_{i = 1, \dots, N}$ are $\mathcal{F}_0$-measurable random $\mathcal{C}(\R_+)$ functions. The model \eqref{eq:Nd-hawkes-general} is a mean-field system of interacting nonlinear Hawkes processes with non-vanishing self-interaction. In the context of neuroscience, \eqref{eq:Nd-hawkes-general} can be seen as a mean-field network of Generalized Linear Model/Spike Response Model neurons. 

Let us now assume that $\mathfrak{h}$ is an Erlang kernel, i.e. there exists $d\in\mathbb{N}^*$ such that
\begin{linenomath}\begin{equation*}
    \mathfrak{h}(t) = c e^{-\alpha t}\frac{t^{d-1}}{(d-1)!}, 
\end{equation*}\end{linenomath}
for some $c \in \R$ and $\alpha > 0$. Self-interaction can then be translated into a `Markovian cascade of memory terms' \cite{DitLoe17, DuaLoe19}: adjusting the initial conditions, \eqref{eq:Nd-hawkes-general} can be equivalently written: for $i = 1, \dots, N$,
\begin{linenomath}\begin{subequations}\label{eq:cascade}
\begin{align}
    Z^{i,N}_t &= \int_{[0,t]\times \R_+} \mathbbm{1}_{\{z\leq f(\M^{i,N}_{s-}(1)+X^{i,N}_{s-})\}}\pi^i(ds,dz), \\
    \M^{i,N}_t &= \M^i_0 + \int_0^t \mathfrak{A} \M^{i,N}_s ds + \mathbf{c} Z^{i,N}_t, \label{eq:M-markov}\\
    X^{i,N}_t &=  \frac{1}{N}\sum_{j\neq i} H^j_t + \frac{1}{N}\sum_{j\neq i}\int_0^t h(t-s)dZ^{j,N}_s,
\end{align}
\end{subequations}\end{linenomath}
where $\M^{i,N}_{s-}(1)$ denotes the first element of the vector $\M^{i,N}_{s-}$. The $d$-by-$d$ matrix $\mathfrak{A}$ has all diagonal terms equal to $-\alpha$, all superdiagonal terms equal to $1$, and all other terms equal to 0; the $d$ dimensional vector $\mathbf{c}$ is defined by $\mathbf{c}(k) = \mathbbm{1}_{k=d}c$. This is equivalent to setting $b(\m) = \mathfrak{A} \m$ and $\Gamma(\m) = (0,\dots, 0, c)$ in \eqref{eq:ALM-Model-M}. The fact that there are $N$ distinct variables $X^{i,N}_t$ instead of a common shared variable $X^N_t$ as in \eqref{eq:ALM-Model} does not affect the mean-field limit since the difference between $X^N_t$ and $X^{i,N}_t$ is of order $1/N$ (see \cite{Che17}).

The generalization of \eqref{eq:cascade} to the case where $\mathfrak{h}$ is a sum of Erlang kernels is straightforward. Of course, a sum of Erlang kernels can simply be a sum of exponential kernels, which is more common in neuroscience \cite{MulBue07,ToyRad09,NauGer12,GerKis14}. Notably, taking one exponential kernel with $c<0$ is enough to model the effects of neuronal refractoriness and spike-frequency adaptation \cite{MulBue07}. 

In the example \eqref{eq:cascade}, the leaky memory variables $\M$ influence the intensity function $f$ but does not influence the interaction function $h$. However, in the general model \eqref{eq:ALM-Model}, $h$ can depend on $\M$. In the context of neuronal modeling, this dependence can be used to account for the effects of short-term synaptic plasticity (STP) \cite{ZucReg02}. Using the notation of the general model \eqref{eq:ALM-Model}, we can describe a network of spiking neurons with refractoriness and `Tsodyks-Markram' STP \cite{TsoPaw98}. The Tsodyks-Markram model \cite{TsoPaw98} captures the interplay between synaptic depression and facilitation and has been used to model working memory \cite{MonBar08}, chaotic dynamics \cite{CorDes13} and learning in hierarchical circuits \cite{PayGue21}. Taking $d=2$, the leaky memory variables (the STP variables of the Tsodyks-Markram model) follow, for initial conditions $\M^i_0$ supported in $[U,1]\times[0,1]$ with $U\in]0,1[$,
\begin{subequations}\label{eq:STP}
\begin{linenomath}\begin{equation}
    \M^{i,N}_t = \M^i_0 + \int_0^t b(\M^{i,N}_s)ds + \int_0^t \Gamma(\M^{i,N}_{s-})dZ^{i,N}_s,
\end{equation}\end{linenomath}
with the vector field
\begin{linenomath}\begin{equation}
    b\left(\m(1),\m(2)\right) = \left(\frac{U-\m(1)}{\tau_F}, \frac{1-\m(2)}{\tau_D}\right),
\end{equation}\end{linenomath}
where $\tau_F>0$ and $\tau_D>0$ are the facilitation and depression timescales respectively, and with the jump function
\begin{linenomath}\begin{equation}\label{eq:STP_Gamma}
    \Gamma\left(\m(1),\m(2)\right) = \left(U(1-\m(1)), -\m(1)\m(2)\right).
\end{equation}\end{linenomath}
\end{subequations}
It is easy to verify that the leaky memory variables $\M^{i,N}_t$ then take values in $[U,1]\times[0,1]$. Finally, we take $f$ independent of the leaky memory variables and $h$ of the form $h(t,a,\m) = \m(1)\m(2)\bar{h}(t)$. The model we just described generalizes interacting Age-dependent Hawkes processes \cite{Che17,RaaDit20} and is more detailed than the model with purely facilitating synapses and without refractoriness studied in \cite{GalLoe20}.

These two motivating examples are clearly special cases of the general model \eqref{eq:ALM-Model}. The fact that in both examples, the variables $\M^{i,N}_t$ relax to some fixed value in the absence of jumps motivates the name `leaky memory'. Importantly, both examples satisfy the main assumptions we will use in this work (see Section~\ref{sec:assumptions}).

\subsection{Methods and relation to previous work}
 To prove propagation of chaos in the mean-field limit, we use the method of coupling \textit{à la Sznitman} \cite{Szn91}: to show the convergence of the time-marginals, we follow Fournier and Löcherbach \cite{FouLoe16} (see also \cite{GalLoe20}); to show the convergence of the processes, we use the method of Delattre, Fournier and Hoffmann \cite{DelFou16} (later used by Chevallier \cite{Che17} and Ditlevsen and Löcherbach \cite{DitLoe17}). Our approach for relating the limit process with the limit PDE differs from previous work \cite{Che17} for it relies on a path integral representation. This representation turns an earlier heuristic method from Naud and Gerstner \cite{NauGer12} into a rigorous argument. Contrarily to \cite{Che17} where PDE solutions in measure space are considered, our method treats PDE solutions in $L^1$ space and does not involve the semigroup theory results of \cite{CanCar13}. More importantly, the path integral method allows us to derive a representation formula for the solution to the PDE. The limit PDE we obtain is a generalization of the Time Elapsed Neuron Network Model of Pakdaman, Perthame and Salort \cite{PakPer10} and of the refractory density equation \cite{GerHem92, Ger95, GerKis14} to the case of neurons with adaptation and short-term synaptic plasticity.
 
\subsection{Plan of the paper}
The main results of this work, namely propagation of chaos (\textbf{Theorem~\ref{theorem:convergence_rate}}) and the characterization of the mean-field limit by a PDE (\textbf{Theorem~\ref{theorem:PDE}}), are presented in Section~\ref{sec:assumptions}, together with the assumptions required. The proof of \textbf{Theorem~\ref{theorem:convergence_rate}} is presented in Section~\ref{sec:theorem_1}. In Section~\ref{sec:alternative}, we show that under more restrictive assumptions, we can get a propagation of chaos result analogous to that of \cite{DelFou16,Che17,DitLoe17}. Finally, the proof of \textbf{Theorem~\ref{theorem:PDE}} is presented in Section~\ref{sec:theorem_2}.

\section{Assumptions and main results}\label{sec:assumptions}

\begin{notations*}
The uniform and Euclidean norms are denoted by $\norm{\cdot}_\infty$ and $\norm{\cdot}$ respectively. We write $\norm{\Gamma}_\infty := \sup_{\m\in\Rd}\norm{\Gamma(\m)}$. We use $C, C_T$ and $C_{T,0}$ to denote positive constants (that can change from line to line) where the subscript $T$ signals the dependence on time and the subscript $0$ the dependence on the law of $\left(A_0^1, \M_0^1, (H^1_t)_{t\in\R_+}\right)$.
\end{notations*}

In this work, we always assume that the functions $f$, $h$, $b$ and $\Gamma$ satisfy:
\begin{assumption}\label{assumption:functions}
\begin{enumerate}[label=(\roman*)]
\item The functions $f$, $h$ and $\Gamma$ are bounded.
\item There exists a bounded, strictly increasing and continuously differentiable function $\psi:\R_+\to\R_+$ with $\psi(0) = 0$ and satisfying 
\begin{linenomath}\begin{equation*}
    |\psi'(a) - \psi'(a^*)|\leq \kappa |\psi(a) - \psi(a^*)|, \qquad \forall a, a^*\in\R_+,
\end{equation*}\end{linenomath}
for some $\kappa>0$, such that, for all $(a, \m, x, a^*, \m^*, x^*) \in (\R_+ \times \R^d \times \R)^2$, and for all $t\in\R_+$,
\begin{linenomath}\begin{align*}
    |f(a, \m, x) - f(a^*, \m^*, x^*_{})| &\leq L_f(|\psi(a) - \psi(a^*)| + \norm{\m - \m^*} + |x - x^*|),  \\
    |h(t, a, \mathbf{m}) - h(t, a^*, \m^*)| &\leq L_h (|\psi(a)-\psi(a^*)| + \norm{\m-\m^*}),
\end{align*}\end{linenomath}
for some $L_f$ and $L_h>0$.
\item The vector field $b$ and jump function $\Gamma$ are Lipschitz continuous.
\end{enumerate}
\end{assumption}

The fact that $f$ is bounded guarantees the well-posedness of the system~\eqref{eq:ALM-Model} and a path-wise unique càdlàg strong solution to~\eqref{eq:ALM-Model} can be constructed using a standard thinning procedure. $(ii)$ says that $f$ and $h$ are Lipschitz continuous with respect to a $\psi$-modified metric on the \textit{age} variable. An example of possible function is $\psi(a) = 1 - \exp(- a \kappa)$. The $\psi$-modified Lipschitz continuity of $f$ implies that the effect of $a$ on $f$ saturates for large $a$. Note that in the STP example \eqref{eq:STP}, since the leaky memory variables take values in the compact $[U,1]\times[0,1]$, the jump function \eqref{eq:STP_Gamma} is effectively Lipschitz.

To prove propagation of chaos, we need some assumptions on $\{(A_0(i), \M_0(i), (H_t(i))_{t\in\R_+})\}_{i \in \mathbb{N}}$:

\begin{assumption} \label{assumption:initial_conditions}
\begin{enumerate}[label=(\roman*)]
\item The 3-tuples $\{(A_0(i), \M_0(i), (H_t(i))_{t\in\R_+})\}_{i \in \mathbb{N}}$ are $i.i.d.$.
\item The random function $(H^1_t)_{t\in\R_+}$ is such that $(\E[H^1_t])_{t\in\R_+} \in \mathcal{C}(\mathbb{R}_+)$.
\item For all $T > 0$, there exists $C_{T,0}>0$ such that $\sup_{t\in[0,T]}\mathrm{Var}[H^1_t] \leq C_{T,0}$.
\end{enumerate}
\end{assumption}
A condition similar to $(iii)$ is used in \cite{Che17}. Note that in \cite{Che17}, Chevallier considers $i.i.d.$ random interaction functions instead of a deterministic function $h$, common to all the point processes. As he proved that, under some square integrability condition, the randomness in the interaction functions averages out in the mean-field limit, we focus here on the fixed $h$ case.

The first main result of this work is a quantified propagation of chaos theorem:
\begin{theorem} \label{theorem:convergence_rate}
Grant Assumptions~\ref{assumption:functions} and~\ref{assumption:initial_conditions}. For all $T>0$, there exists $C_{T,0}>0$ such that
\begin{linenomath}\begin{equation} \label{eq:convergence_A_M}
    \sup_{t \in [0,T]}\mathbb{E}\left[|\psi(A^{1,N}_t) - \psi(A^1_t)|+\|\M^{1,N}_t - \M^1_t\| + |X^N_t - x_t|\right] \leq \frac{C_{T,0}}{\sqrt{N}},
\end{equation}\end{linenomath}
where $(A^1_t, \M^1_t, x_t)_{t\in\R_+}$ (the limit process) is given by the path-wise unique càdlàg strong solution to
\begin{subequations}\label{eq:limit}
\begin{linenomath}\begin{align}
    A^1_t &= A^1_0 + t - \int_0^t A^1_{s-}dZ^1_s, \label{eq:limit_A}\\
    \M^1_t &= \M^1_0 + \int_0^t b(\M^1_s)ds + \int_0^t \Gamma(\M^1_{s-})dZ^1_s,\label{eq:limit_M}\\
    Z^1_t &= \int_{[0,t]\times\R_+}\mathbbm{1}_{z\leq f(A^1_{s-}, \M^1_{s-}, x_s)}\pi^1(ds,dz), \label{eq:limit_Z}\\
    x_t &= \E[H^1_t] + \int_0^t \E[h(t-s,A^1_s,\M^1_s)f(A^1_s,\M^1_s,x_s)]ds.
\end{align}\end{linenomath}
\end{subequations}
Furthermore, for all $t\in[0,T]$, writing $\mathcal{L}\left(\psi(A^1_t), \M^1_t\right)$ the law of $(\psi(A^1_t), \M^1_t)$, there exists $C'_{T,0}>0$ such that
\begin{linenomath}\begin{equation}\label{eq:empirical}
    \sup_{t \in [0,T]}\E\left[W_1\left(\frac{1}{N}\sum_{i=1}^N\delta_{(\psi(A^{i,N}_t),\M^{i,N}_t)}, \mathcal{L}\left(\psi(A^1_t), \M^1_t\right)\right)\right] \leq \frac{C'_{T,0}}{\sqrt{N}},
\end{equation}\end{linenomath}
where $W_1$ denotes the $1$-Wasserstein distance.
\end{theorem}
It directly follows from Theorem~\ref{theorem:convergence_rate} and the Continuous mapping theorem that for all $t>0$, the empirical measure of the system \eqref{eq:ALM-Model} at time $t$ converges in probability to the time-marginal of the law of the limit process:
\begin{linenomath}\begin{equation}\label{eq:time-marginals}
    \frac{1}{N}\sum_{i=1}^N\delta_{(A^{i,N}_t,\M^{i,N}_t, X^N)} \xrightarrow[N\to\infty]{\mathbb{P}} \mathcal{L}\left(A^1_t, \M^1_t, x_t\right).
\end{equation}\end{linenomath}

The second main result relates the time-marginals $\mathcal{L}\left(A^1_t, \M^1_t, x_t\right)$ of the law of the limit process with the solution to a nonlocal transport equation. To formulate the transport equation, we define the \textit{jump mapping} $\gamma(\m) = \m + \Gamma(\m)$ and we write $\nabla\cdot$ the divergence operator in $\Rd$. To stay within the standard framework of (mass-conservative) transport equation with solutions in $L^1$ \cite{Per07}, we need:
\begin{assumption}\label{assumption:gamma_dif}
\begin{enumerate}[label=(\roman*)]
\item The vector field satisfy $b \in \mathcal{C}^1(\Rd, \Rd)$ and $\nabla\cdot b \in \mathcal{C}^1(\Rd, \Rd)$.
\item The jump mapping $\gamma$ is a proper local $\mathcal{C}^1$-diffeomorphism.
\end{enumerate}
\end{assumption} 

\begin{theorem} \label{theorem:PDE}
Grant Assumptions~\ref{assumption:functions} and \ref{assumption:gamma_dif}. Further assume that the law of the initial condition $(A^1_0,\M^1_0)$ is the absolutely continuous probability measure $u_0(a,\m)dad\m$ and $(\E[H_t])_{t\in\R_+} = (\bar{H}^1_t)_{t\in\R_+}\in \mathcal{C}(\R_+)$. Then, the time-marginals $\rho_t\otimes\delta_{x_t}:=\mathcal{L}\left(A^1_t, \M^1_t, x_t\right)$ of the law of the limit process \eqref{eq:limit} correspond the unique weak solution to
\begin{subequations}\label{eq:PDE}
\begin{linenomath}\begin{align}
    &\partial_t \rho_t(a,\m) + \partial_a \rho_t(a,\m) + \nabla\cdot\left(b(\m)\rho_t(a,\m)\right) = -f(a,\m,x_t)\rho_t(a,\m),\\
    &\rho_t(0,\cdot) = \gamma_*\left(\int_{\R_+} f(a,\cdot,x_t)\rho_t(a,\cdot)da\right), \label{eq:PDE_border}\\
    &x_t = \bar{H}_t + \int_0^t\int_{\Rd}\int_{\R_+}h(t-s,a,\m)f(a,\m,x_s)\rho_s(a,\m)dad\m ds,\\
    &\rho_0(a,\m) = u_0(a,\m),
\end{align}\end{linenomath}
\end{subequations}
(where $\nabla \cdot$ denotes the divergence on the variables $\m$ and $\gamma_*(\dots)$ denotes the pushforward measure by $\gamma$) in the sense that $(\rho,x) \in \mathcal{C}(\mathbb{R}_+, L^1(\mathbb{R}_+ \times \mathbb{R}^d)) \times \mathcal{C}(\mathbb{R}_+)$ and, for all $G \in \mathcal{C}^\infty_c(\mathbb{R}_+ \times \mathbb{R}_+ \times \mathbb{R}^d)$,
\begin{linenomath}\begin{multline} \label{eq:weak}
    0=\int_{\R^d}\int_{\R_+} G(0,a, \m)u_0(a,\m)dad\mathbf{m} + \int_{\R_+}\int_{\R^d}\int_{\R_+}\Big\{[\partial_t + \partial_a + b(\m)\cdot\nabla] G(t,a,\m) \\
    + (G(t,0, \gamma(\m))-G(t,a, \m))f(a,\m,x_t)\Big\}\rho_t(a,\m)dad\m dt,
\end{multline}\end{linenomath}
where $\nabla$ denotes the gradient operator on the variables $\m$.
\end{theorem}

Assumption~\ref{assumption:gamma_dif}~\textit{(ii)} guarantees that for all $\m\in\gamma(\Rd)$, the preimage $\gamma^{-1}(\m)$ is a finite set of points and \eqref{eq:PDE_border} can be more explicitly written
\begin{linenomath}\begin{equation}\label{eq:border_explicit}
    \rho_t(0,\m) = \mathbbm{1}_{\gamma(\Rd)}(\m)\sum_{\m'\in\gamma^{-1}(\m)}\frac{1}{|\det(\mathbf{J}_\gamma(\m'))|}\int_{\R_+}f(a,\m',x_t)\rho_t(a,\m')da,
\end{equation}\end{linenomath}
where $\det(\mathbf{J}_\gamma(\m'))$ denotes the determinant of the Jacobian matrix $\mathbf{J}_\gamma(\m')$.

All the results and proofs can be adapted to the simpler case where the system is not age-dependent, as in the Erlang kernel example~\eqref{eq:cascade}. For this example, the limit PDE~\eqref{eq:PDE} becomes
\begin{subequations}\label{eq:PDE_cascade}
\begin{linenomath}\begin{align*}
    &\partial_t \rho_t(\m) \rho_t(\m) + \mathfrak{A}\nabla\cdot\left(\m\rho_t(\m)\right) = f(\m(1) + x_t)\rho_t(\m-\mathbf{c})-f(\m(1) + x_t)\rho_t(\m),\\
    &x_t = \bar{H}_t + \int_0^t h(t-s)\int_{\Rd}f(\m(1) + x_s)\rho_s(\m)d\m ds,\\
    &\rho_0(\m) = u_0(\m).
\end{align*}\end{linenomath}
\end{subequations}

For the STP example~\eqref{eq:STP}, the limit PDE reads
\begin{subequations}
\begin{linenomath}\begin{align*}
    &\partial_t \rho_t(a,\m) + \partial_a \rho_t(a,\m) + \nabla\cdot\left(b(\m)\rho_t(a,\m)\right) = -f(a,x_t)\rho_t(a,\m),\\
    &\rho_t(0,\m) = \mathbbm{1}_{\gamma(D)}(\m)\frac{1}{1-\m(1)}\int_{\R_+} f(a,x_t)\rho_t(a,\gamma^{-1}(\m))da,\\
    &x_t = \bar{H}_t + \int_0^t\bar{h}(t-s)\int_{\Rd}\int_{\R_+}\m(1)\m(2)f(a,x_s)\rho_s(a,\m)dad\m ds,\\
    &\rho_0(a,\m) = u_0(a,\m),
\end{align*}\end{linenomath}
\end{subequations}
where $D = ]U,1[\times]0,1[$ and $\gamma(\m) = \left(U + [1-U]\m(1), \,[1-\m(1)]\m(2)\right)$.\footnote{Using \eqref{eq:border_explicit}, a simple calculation gives $\det (\mathbf{J}_\gamma(\gamma^{-1}(\m))) = 1-\m(1)$.}

Theorems~\ref{theorem:convergence_rate} and \ref{theorem:PDE} have two important implications for neuronal modeling: first, they provide a rigorous footing to multidimensional population density equations, which could be simulated using mesh-based methods described in \cite{MulBue07,LaiDek17,DekLai19}; second, they confirm that, not only in the simple `renewal' cases, the PDE point of view can be used to study the nonlinear dynamics of large networks of spiking neurons \cite{FonSch21}. 

\section{Proof of Theorem~\ref{theorem:convergence_rate} (Propagation of chaos)}\label{sec:theorem_1}
The approach here is standard. We use a fixed-point argument to show that the limit process~\eqref{eq:limit} is well-defined. Then, we use the coupling method \cite{Szn91} to prove that a \textit{typical particle} converges to the limit process.

\subsection{Well-posedness of the limit process}

\begin{proposition} \label{proposition:well_posedness_ALM-SDE}
Grant Assumption~\ref{assumption:functions} and assume that $(\E[H^1_t])_{t\in\R_+} \in \mathcal{C}(\mathbb{R}_+)$. There exists a path-wise unique càdlàg strong solution $(A^1_t, \M^1_t, x_t)_{t \in \mathbb{R}_+}$ taking values in $\mathbb{R}_+ \times \mathbb{R}^d \times \mathbb{R}$ to the limit equation~\eqref{eq:limit}. Furthermore, $(x_t)_{t\in\R_+} \in \mathcal{C}(\R_+)$.
\end{proposition}
\begin{proof}
For all $y\in \mathcal{C}(\R_+)$, let us write $(A^y_t, \M^y_t, x^y_t)_{t\in\R_+}$ the càdlàg strong solution to
\begin{linenomath}\begin{align*}
    A^y_t &= A^1_0 + t - \int_0^t A^y_{s-}dZ^y_s, \\
    \M^y_t &= \M^1_0 + \int_0^t b(\M^y_s)ds + \int_0^t \Gamma(\M^y_{s-})dZ^y_s,\\
    Z^y_t &= \int_{[0,t]\times\R_+}\mathbbm{1}_{z\leq f(A^y_{s-}, \M^y_{s-}, y_s)}\pi^1(ds,dz).
\end{align*}\end{linenomath}
Then, we set
\begin{linenomath}\begin{equation*}
    x^y_t = \E[H^1_t] + \int_0^t \E[h(t-s,A^y_s,\M^y_s)f(A^y_s,\M^y_s,y_s)]ds.
\end{equation*}\end{linenomath}

Since $f$ and $h$ are bounded, by dominated convergence, we have that $(x^y_t)_{t\in\R_+}\in\mathcal{C}(\R_+)$. Thus, for all $T>0$, we can define the operator
\begin{linenomath}\begin{equation*}
    \Phi_{T} : \mathcal{C}([0, T]) \to \mathcal{C}([0, T]), \qquad (y_t)_{t \in [0, T]} \mapsto (x^y_t)_{t \in [0, T]}.
\end{equation*}\end{linenomath}
The solution $(A^y_t, \M^y_t, x^y_t)_{t\in[0,T]}$ is a solution to \eqref{eq:limit} on $[0,T]$ if and only if $(y_t)_{t \in [0, T]} = (x^y_t)_{t \in [0, T]}$, or equivalently, if and only if $(y_t)_{t \in [0, T]}$ is a fixed point of $\Phi_T$. We are going to show that for $T$ small enough, $\Phi_T$ is a contraction for the uniform norm. 

For all $y$ and $y^*\in\mathcal{C}([0,T])$, by triangular inequality and using the Lipschitz continuity and the boundedness of $f$ and $h$, we have, for all $t\in[0,T]$,
\begin{linenomath}\begin{equation*}
    \left|\int_0^t\mathbb{E}[h(t-s, A_s^y, \M_s^y)f(A_s^y, \M_s^y, y_s)]ds - \int_0^t\mathbb{E}[h(t-s, A_s^{y^*}, \M_s^{y^*})f(A_s^{y^*}, \M_s^{y^*}, y^*_s)]ds\right|\leq C \int_0^t \Delta_s ds,
\end{equation*}\end{linenomath}
where 
\begin{linenomath}\begin{equation*}
    \Delta_s := \E\left[|\psi(A^y_s) - \psi(A^{y^*}_s)| + \|\M^y_s - \M^{y^*}_s\| + |y_s - y^*_s|\right].
\end{equation*}\end{linenomath}
By Itô's formula for jump processes, 
\begin{linenomath}\begin{equation*}
    \psi(A^y_t) = \psi(A^1_0) + \int_0^t\psi'(A_s^y)ds - \int_{[0,t]\times \R_+}\psi(A_{s-}^y)\mathbbm{1}_{z\leq f(A^y_{s-},\M^y_{s-}, y_s)}\pi^1(ds,dz).
\end{equation*}\end{linenomath}
Whence,
\begin{linenomath}\begin{multline*}
    \mathbb{E}[|\psi(A^y_t) - \psi(A^{y^*}_t)|] \leq \mathbb{E}\left[\left|\int_0^t \psi'(A^y_s) - \psi'(A^{y^*}_s)ds\right|\right] \\
     + \mathbb{E}\left[\left|\int_{[0,t]\times\R_+}\psi(A^y_{s-}) \mathbbm{1}_{z \leq f(A^y_{s-}, \M^y_{s-}, y_s)} - \psi(A^{y^*}_{s-}) \mathbbm{1}_{z \leq f(A^{y^*}_{s-}, \M^{y^*}_{s-}, y^*_s)}\pi^1(ds,dz)\right|\right].
\end{multline*}\end{linenomath}
Notice that by Assumption~\ref{assumption:functions}, $|\psi'(A^y_s) - \psi'(A^{y^*}_s)|\leq \kappa |\psi(A^y_s) - \psi(A^{y^*}_s)|$. Then, by triangular inequality and using the Lipschitz continuity and the boundedness of $f$ and $\psi$, we easily get
\begin{linenomath}\begin{equation*}
    \mathbb{E}[|\psi(A^y_t) - \psi(A^{y^*}_t)|] \leq C \int_0^t \Delta_s ds;
\end{equation*}\end{linenomath}
similarly, using the Lipschitz continuity of $b$, $\Gamma$ and $f$ and the the boundedness of $\Gamma$ and $f$, we get $\mathbb{E}[\|\M^y_t - \M^{y^*}_t\|] \leq C \int_0^t \Delta_s ds$.
Thus, for all $t\in[0,T]$,
\begin{linenomath}\begin{equation*}
    \Delta_t \leq C\int_0^t\Delta_s ds + \norm{y - y^*}_\infty,
\end{equation*}\end{linenomath}
and by Grönwall's lemma,
\begin{linenomath}\begin{equation*}
    \Delta_t \leq \norm{y - y^*}_\infty\exp(Ct).
\end{equation*}\end{linenomath}
Whence,
\begin{linenomath}\begin{equation*}
    \norm{\Phi_T(y)-\Phi_T(y^*)}_\infty \leq C'T\exp(CT)\norm{y-y^*}_\infty.
\end{equation*}\end{linenomath}
For $T$ small enough, $\Phi_T$ is a contraction and has a unique fixed point by Banach's fixed-point theorem. The fixed point gives the unique solution to \eqref{eq:limit} on $[0,T]$. Since the constants $C$ and $C'$ do not depend on $T$ nor on the law of $(A^1_0,\M^1_0,(H^1_t)_{t\in\R_+})$, we can iterate the argument above on successive time intervals of length $T$ to obtain the solution to~\eqref{eq:limit} on $\R_+$.
\end{proof}

\subsection{Convergence}\label{sec:convergence}
\begin{proof}
The proof of the convergence \eqref{eq:convergence_A_M} follows the same general strategy as in \cite[Theorem~7]{DelFou16}.

For all $i=1,\dots,N$, we define the \textit{coupled} limit process process $(A^i_t, \M^i_t, Z^i_t)_{t\in\R_+}$ as the path-wise unique càdlàg strong solution to  
\begin{linenomath}\begin{align*}
    A^i_t &= A^i_0 + t - \int_0^t A^i_{s-}dZ^i_s, \\
    \M^i_t &= \M^i_0 + \int_0^t b(\M^i_s)ds + \int_0^t \Gamma(\M^i_{s-})dZ^i_s,\\
    Z^i_t &= \int_{[0,t]\times\R_+}\mathbbm{1}_{z\leq f(A^i_{s-}, \M^i_{s-}, x_s)}\pi^i(ds,dz), \\
    x_t &= \E[H^i_t] + \int_0^t \E[h(t-s,A^i_s,\M^i_s)f(A^i_s,\M^i_s,x_s)]ds.
\end{align*}\end{linenomath}
The process $(A^i_t, \M^i_t, Z^i_t)_{t\in\R_+}$ is coupled to $(A^{i,N}_t, \M^{i,N}_t, Z^{i,N}_t)_{t\in\R_+}$ in the sense that it shares the same $(A^i_0, \M^i_0, (H^i_t)_{t\in\R_+})$ and the same Poisson random measure $\pi^i$. The variable $x_t$ has no index $i$ as it the same for all $i$; it can be interpreted as the deterministic time-varying `mean field' which acts uniformly on all the individual processes. Importantly, the limit processes $\{(A^i_t, \M^i_t, Z^i_t)_{t\in\R_+}\}_{i=1}^N$ are \textit{i.i.d.}.

For all $t\geq 0$, let us define
\begin{linenomath}\begin{equation*}
    \Delta^{1,N}_t := \E\left[|\psi(A^{1,N}_t) - \psi(A^1_t)| + \|\M^{1,N}_t - \M^1_t\| + |X^N_t - x_t|\right]. 
\end{equation*}\end{linenomath}
Arguing as in the proof of Proposition~\ref{proposition:well_posedness_ALM-SDE}, we get
\begin{linenomath}\begin{equation*}
    \E\left[|\psi(A^{1,N}_t) - \psi(A^1_t)| + \|\M^{1,N}_t - \M^1_t\|\right] \leq C \int_0^t \Delta^{1,N}_s ds.
\end{equation*}\end{linenomath}
It remains to control the term $\E\left[|X^N_t - x_t|\right]$:

Fix $T>0$. For all $t\in[0,T]$, by triangular inequality,
\begin{linenomath}\begin{align}
    &\E\left[|X^N_t - x_t|\right] \leq \mathbb{E}\biggr[\biggr|\frac{1}{N}\sum_{i=1}^N  H^i_t - \mathbb{E}[H^1_t]\biggr|\biggr] \nonumber\\
    &+\mathbb{E}\biggr[\biggr| \frac{1}{N}\sum_{i=1}^N \int_{[0,t]\times \R_+} h(t-s, A_{s-}^{i,N}, \M_{s-}^{i,N})\mathbbm{1}_{z \leq f(A^{i,N}_{s-}, \M_{s-}^{i,N}, X^N_{s-})}\pi^i(ds,dz) \nonumber\\
    &\qquad \qquad \qquad \qquad \qquad  - \frac{1}{N}\sum_{i=1}^N\int_{[0,t]\times \R_+}  h(t-s, A^i_{s-}, \M^i_{s-})\mathbbm{1}_{z \leq f(A^i_{s-}, \M^i_{s-}, x_s)}\pi^i(ds,dz)\biggr|\biggr] \nonumber\\
    &+\mathbb{E}\biggr[\biggr| \frac{1}{N}\sum_{i=1}^N \int_{[0,t]\times \R_+} h(t-s, A^i_{s-}, \M^i_{s-})\mathbbm{1}_{z \leq f(A^i_{s-}, \M^i_{s-}, x_s)}\pi^i(ds,dz) \nonumber\\
    &\qquad \qquad \qquad \qquad \qquad \qquad \qquad \qquad \qquad  \quad - \int_0^t\mathbb{E}[h(t-s,A^1_s,\M^1_s)f(A^1_s, \M^1_s, x_s)]ds\biggr|\biggr] \nonumber\\
    =:& Q^N_t + R^N_t + S^N_t. \label{eq:bound_QRS}
\end{align}\end{linenomath}
By Cauchy-Schwarz inequality and Assumption~\ref{assumption:initial_conditions}, 
\begin{linenomath}\begin{equation*}
    Q^N_t \leq \left(\frac{\mathrm{Var}[H^1_t]}{N}\right)^\frac{1}{2} \leq \frac{C_{T,0}}{\sqrt{N}}.
\end{equation*}\end{linenomath}
By exchangeability, triangular inequality and by the Lipschitz continuity and boundedness of $f$ and $h$,
\begin{linenomath}\begin{align*}
    R^N_t \leq& C\int_0^t \Delta^{1,N}_s ds.
\end{align*}\end{linenomath}
By Cauchy-Schwarz inequality,
\begin{linenomath}\begin{align*}
    S^N_t \leq& \mathbb{E}\biggr[\biggr( \frac{1}{N}\sum_{i=1}^N \int_{[0,t]\times \R_+} h(t-s, A^i_{s-}, \M^i_{s-})\mathbbm{1}_{z \leq f(A^i_{s-}, \M^i_{s-}, x_s)}\pi^i(ds,dz) \\
    &\qquad \qquad \qquad \qquad \qquad \qquad \qquad \qquad \qquad \quad - \int_0^t\mathbb{E}[h(t-s,A^1_s,\M^1_s)f(A^1_s, \M^1_s, x_s)]ds\biggr)^2\biggr]^\frac{1}{2} \\
    =&\frac{1}{\sqrt{N}}\mathrm{Var}\left[\int_{[0,t]\times\R_+} h(t-s,A^1_{s-},\M^1_{s-})\mathbbm{1}_{z \leq f(A^1_{s-}, \M^1_{s-}, x_s)}\pi^1(ds,dz)\right]^{\frac{1}{2}}.
\end{align*}\end{linenomath}

However, writing $\tilde{\pi}^1(ds,dz):= \pi^1(ds,dz) - dsdz$ the compensated Poisson random measure, we have, by Itô isometry for compensated jump processes,
\begin{linenomath}\begin{multline*}
    \mathrm{Var}\left[\int_{[0,t]\times\R_+} h(t-s,A^1_{s-},\M^1_{s-})\mathbbm{1}_{z \leq f(A^1_{s-}, \M^1_{s-}, x_s)}\pi^1(ds,dz)\right] \\
    = \E\left[\left(\int_{[0,t]\times\R_+} h(t-s,A^1_{s-},\M^1_{s-})\mathbbm{1}_{z \leq f(A^1_{s-}, \M^1_{s-}, x_s)}\tilde{\pi}^1(ds,dz)\right)^2\right] \leq T \norm{h}_\infty^2 \norm{f}_\infty.
\end{multline*}\end{linenomath}
Hence, $S^N_t \leq N^{\frac{1}{2}}C_T$. Gathering the bounds, we get
\begin{linenomath}\begin{equation*}
    \E\left[|X^N_t - x_t|\right] \leq C\int_0^t\Delta^{1,N}_s ds + \frac{C_{T,0}}{\sqrt{N}}.
\end{equation*}\end{linenomath}
Finally,
\begin{linenomath}\begin{equation*}
    \Delta^{1,N}_t \leq C\int_0^t\Delta^{1,N}_s ds + \frac{C_{T,0}}{\sqrt{N}}, \qquad \forall t \in [0, T],
\end{equation*}\end{linenomath}
and by Gr\"{o}nwall's lemma, 
\begin{linenomath}\begin{equation}\label{eq:delta_conv}
    \Delta^{1,N}_t \leq \frac{C_{T,0}}{\sqrt{N}}\exp(C_T T), \qquad \forall t \in [0,T],
\end{equation}\end{linenomath}
which concludes the proof of \eqref{eq:convergence_A_M}.

By echangeability,
\begin{linenomath}\begin{multline*}
    \E\left[W_1\left(\frac{1}{N}\sum_{i=1}^N\delta_{(\psi(A^{i,N}_t),\M^{i,N}_t)}, \mathcal{L}\left(\psi(A^1_t), \M^1_t\right)\right)\right] \\
    \leq \E\left[|\psi(A^{1,N}_t - \psi(A^1_t)| + \|\M^{1,N}_t - \M^1_t\|\right] 
    + \E\left[W_1\left(\frac{1}{N}\sum_{i=1}^N\delta_{(\psi(A^{i}_t),\M^{i}_t)}, \mathcal{L}\left(\psi(A^1_t), \M^1_t\right)\right)\right].
\end{multline*}\end{linenomath}
Then, we simply use \eqref{eq:delta_conv} and a result on the convergence of the empirical measures in Wasserstein distance \cite[Theorem~1]{FouGui15} to get \eqref{eq:empirical}.
\end{proof}

\section{Alternative propagation of chaos result}\label{sec:alternative}
Theorem~\ref{theorem:convergence_rate} guarantees the convergence of the time-marginals (see \eqref{eq:time-marginals}), which is sufficient for relating the empirical measure of the system~\eqref{eq:ALM-Model} with the PDE \eqref{eq:PDE}. However, under more restrictive assumptions on the vector field $b$ and the jump mapping $\gamma$, it is possible to get the convergence of the processes, as in \cite{DelFou16,Che17,DitLoe17}. 

\begin{assumption}\label{assumption:alternative}
    \begin{enumerate}[label=(\roman*)]
    \item Writing $(B_t)_{t\in\R_+}$ the flow associated with the vector field $b$, for all $t\geq 0$, $B_t$ is $1$-Lipschitz for the Euclidean distance.
    \item The jump mapping $\gamma$ is $1$-Lipschitz for the Euclidean distance.
    \end{enumerate}
\end{assumption}
\begin{theorem} \label{theorem:alternative}
Grant Assumptions~\ref{assumption:functions},~\ref{assumption:initial_conditions} and \ref{assumption:alternative}. For all $T>0$, there exists $C_{T,0}>0$ such that
\begin{linenomath}\begin{equation} \label{eq:convergence_A_M_alternative}
    \mathbb{E}\left[\sup_{t \in [0,T]}|\psi(A^{1,N}_t) - \psi(A^1_t)|+\|\M^{1,N}_t - \M^1_t\|\right] \leq \frac{C_{T,0}}{\sqrt{N}},
\end{equation}\end{linenomath}
where $(A^1_t, \M^1_t)_{t\in\R_+}$ is given by the path-wise unique strong solution to \eqref{eq:limit}.
\end{theorem}
\begin{proof}
The well-posedness of the limit process $(A^1_t, \M^1_t, x_t)_{t\in\R_+}$ has already been proved in Proposition~\ref{proposition:well_posedness_ALM-SDE}. For the convergence, we follow the same strategy as in \cite[Theorem~8]{DelFou16} (see also \cite[Theorem~IV.1]{Che17} and \cite[Theorem~1]{DitLoe17}).

Let $\{(A^i_t, \M^i_t, Z^i_t)_{t\in\R_+}\}_{i=1}^N$ be the same coupled limit process as in Section~\ref{sec:convergence}. The integral $\int_0^t |d(Z^{1,N}_s - Z^1_s)|$ counts the number of times one counting process jumps whereas the other does not, on the time interval $[0,t]$. We define
\begin{linenomath}\begin{equation*}
    \delta_t^N := \mathbb{E}\left[\int_0^t |d(Z^{1,N}_s - Z^1_s)|\right] = \int_0^t \mathbb{E}[|f(A^{1,N}_s, \M^{1,N}_s, X^N_s) - f(A^1_s, \M^1_s, x_s)|]ds.
\end{equation*}\end{linenomath}
The key observation is that Assumption~\ref{assumption:alternative} guarantees
\begin{linenomath}\begin{equation} \label{eq:bound_A_M_delta}
    \quad \mathbb{E}\left[\sup_{s \in [0,t]}|\psi(A^{1,N}_s) - \psi(A^1_s)|+\|\M^{1,N}_s - \M^1_s\|_2\right] \leq C\delta_t^N:
\end{equation}\end{linenomath}
Clearly, $\sup_{s \in [0,t]}|\psi(A^{1,N}_s) - \psi(A^1_s)| \leq \norm{\psi}_\infty \cdot \mathbbm{1}_{\int_0^t |d(Z^{1,N}_s - Z^1_s)|>0} \leq \norm{\psi}_\infty\int_0^t |d(Z^{1,N}_s - Z^1_s)|$, which implies that $\mathbb{E}[\sup_{s \in [0,t]}|\psi(A^{1,N}_s) - \psi(A^1_s)|] \leq \norm{\psi}_\infty \delta_t^N$. On the other hand, by Assumption~\ref{assumption:alternative}~\textit{(i)}, in a time interval with no jumps in $Z^{1,N}$ nor $Z^1$, $\|\M_t^{1,N} - \M^1_t\|_2$ can not increase. If both $Z^{1,N}$ and $Z^1$ jump at time $t$, by Assumption~\ref{assumption:alternative}~\textit{(ii)}, $\|\M_t^{1,N} - \M^1_t\| \leq \|\M_{t-}^{1,N} - \M^1_{t-}\|$. Hence, the only way to have $\|\M_t^{1,N} - \M^1_t\|_2>\|\M_{t-}^{1,N} - \M^1_{t-}\|_2$ is if $Z^{1,N}$ jumps at time $t$ but not $Z^1$ or vice versa. However, in these cases, the increase is bounded by $\norm{\Gamma}_\infty$. In summary, we have that $\sup_{s \in [0,t]}\|\M^{1,N}_s - \M^1_s\|_2 \leq \norm{\Gamma}_\infty\int_0^t |d(Z^{1,N}_s - Z^1_s)|$, which concludes the verification of \eqref{eq:bound_A_M_delta}.

We now have to control $\delta_t^N$. Using~\eqref{eq:bound_A_M_delta}, we have that
\begin{linenomath}\begin{multline*}
    \delta_t^N \leq L_f \int_0^t \mathbb{E}\left[|\psi(A^{1,N}_s) - \psi(A^1_s)| + \|\M^{1,N}_s - \M^1_s\|_2 + |X^N_s - x_s|\right]ds \\
    \leq L_f \int_0^t C \delta_s^N + \mathbb{E}[|X^N_s - x_s|]ds.
\end{multline*}\end{linenomath} 
Fix $T>0$, for all $s\in[0,T]$, we can bound $\mathbb{E}[|X^N_s - x_s|]$ as in the proof of Theorem~\ref{theorem:convergence_rate} (see \eqref{eq:bound_QRS}):
\begin{linenomath}\begin{equation*}
    \mathbb{E}[|X^N_s - x_s|] \leq Q^N_s + R^N_s + S^N_s.
\end{equation*}\end{linenomath} 
with the same variance bounds $Q^N_s + S^N_s \leq N^{\frac{1}{2}}C_{T,0}$. By exchangeability, triangular inequality and using \eqref{eq:bound_A_M_delta}, we have
\begin{linenomath}\begin{equation*}
    R^N_s \leq \norm{h}_\infty \delta_s^N +\; \norm{f}_\infty L_h \int_0^s \mathbb{E}[|\psi(A_u^{1,N}) - \psi(A^1_u)| - \|\M_u^{1,N} - \M^1_u\|]du \leq C_T \delta_s^N.
\end{equation*}\end{linenomath} 
Gathering the bounds, we get
\begin{linenomath}\begin{equation*}
    \delta_t^N \leq C_T \int_0^t \delta_s^N ds + \frac{C_{T,0}}{\sqrt{N}}, \qquad \forall t \in [0, T],
\end{equation*}\end{linenomath}
and we conclude using Gr\"{o}nwall's lemma.
\end{proof}

By exchangeability, \eqref{eq:convergence_A_M_alternative} implies that for all $T>0$, there exists $C'_{T,0}>0$ such that
\begin{linenomath}\begin{equation*} 
    \mathbb{E}\left[\sup_{t \in [0,T]}|\psi(A^{1,N}_t) - \psi(A^1_t)|+\|\M^{1,N}_t - \M^1_t\|_2 + |\psi(A^{2,N}_t) - \psi(A^2_t)|+\|\M^{2,N}_t - \M^2_t\|_2\right] \leq \frac{C'_{T,0}}{\sqrt{N}}.
\end{equation*}\end{linenomath}
By standard arguments on the Skorokhod metric and the Continuous mapping theorem, we have the weak convergence
\begin{linenomath}\begin{equation*}
    \left(A^{1,N}_t, \M^{1,N}_t, A^{2,N}_t, \M^{2,N}_t\right)_{t\in\R_+} \xrightarrow[N\to\infty]{w}\left(A^1_t, \M^1_t, A^2_t, \M^2_t\right)_{t\in\R_+}.
\end{equation*}\end{linenomath}
Since, $(A^1_t, \M^1_t)_{t\in\R_+}$ and $(A^2_t, \M^2_t)_{t\in\R_+}$ have the same law, by \cite[Proposition~2.2]{Szn91}, we have the convergence in probability of the empirical measure of the system \eqref{eq:ALM-Model} to the law of the limit process:
\begin{linenomath}\begin{equation} \label{eq:empirical_alternative}
\frac{1}{N}\sum_{i=1}^N \delta_{(A^{i,N}_t, \M^{i,N}_t)_{t \in \R_+}} \xrightarrow[N \to \infty]{\mathbb{P}} \mathcal{L}\left((A^1_t, \M^1_t)_{t \in \R_+}\right) \quad\text{ in } \mathcal{P}(\mathcal{D}(\R_+, \R_+ \times \Rd)),
\end{equation}\end{linenomath}
where $\mathcal{P}(\mathcal{D}(\R_+, \R_+ \times \Rd))$ denotes the space of probability measures on the Skorokhod space $\mathcal{D}(\R_+, \R_+ \times \Rd)$ of càdlàg functions $\R_+ \to \R_+ \times \Rd$ and $\mathcal{L}\left(\left(A^1_t, \M^1_t\right)_{t \in \R_+}\right)$ denotes the law of the process $\left(A^1_t, \M^1_t\right)_{t \in \R_+}$. 

The convergence \eqref{eq:empirical_alternative} is clearly stronger than the convergence of the time-marginals \eqref{eq:empirical} but it requires the additional Assumption~\ref{assumption:alternative}, which is somewhat restrictive.

\section{Proof of Theorem~\ref{theorem:PDE} (Transport equation for the empirical measure)}\label{sec:theorem_2}
Here, our aim is to show that if we write $(\rho_t \otimes \delta_{x_t})_{t\in\R_+}$ the time-marginals of the law of the process~\eqref{eq:limit}, then $(\rho_t,x_t)_{t\in\R_+}$ is a weak solution to~\eqref{eq:PDE}. To show this, we use the limit process to construct a representation formula for $\rho_t$. The representation formula is obtained by making rigorous the heuristic `path integral' method in \cite{NauGer12}. We then show that the path integral representation gives a weak solution to~\eqref{eq:PDE}. Finally, we prove that weak solution to~\eqref{eq:PDE} is unique. 

\subsection{Path integral representation for the time-marginals of the law ot the limit process}
To formulate the path integral representation, we first need to introduce some notations and definitions.

Let $(x_t)_{t\in\R_+}$ be given by the limit process~\eqref{eq:limit} and let us write $(A^*_t,\M^*_t)_{t\in\R_+}$ the càdlàg process following the dynamics~\eqref{eq:limit_A} and \eqref{eq:limit_M} given $(x_t)_{t\in\R_+}$ and the initial condition $(a_0, \m_0)\in \R_+\times\Rd$. For all $t>0$, $(A^*_t,\M^*_t)$ is deterministic given the initial condition and the jump times in $[0,t]$. Hence, for all $k\in\mathbb{N}$ (the number of jumps in $[0,t]$) and for all $0<t_1<\dots<t_k\leq t$ (the jump times in $[0,t]$), we can define, recursively, the mappings $\theta^k_t(t_1,\dots,t_k): \Rd \to \Rd$, giving $\M^*_t$ as a function of the initial condition $\m_0$:
\begin{subequations} \label{eq:def_theta}
\begin{linenomath}\begin{align}
    \theta^0_t &:= B_t, \\
    \forall k\geq 1,\quad \theta^k_t(t_1, \dots, t_k) &:= B_{t-t_k}\circ\gamma\circ \theta^{k-1}_{t_k}(t_1, \dots, t_{k-1}),
\end{align}\end{linenomath}
\end{subequations}
where $(B_t)_{t\in\R_+}$ is the flow associated with the vector field $b$.

For all $k\geq 1$, we can now define the mapping
\begin{linenomath}\begin{equation*}
    \phi^k_t\left(\begin{pmatrix} t_1 \\ \vdots \\ t_{k-1} \\ t_k \\ \m_0 \end{pmatrix}\right) =
    \begin{pmatrix} t_1 \\ \vdots \\ t_{k-1} \\ t-t_k \\ \theta^k_t(t_1, \dots, t_{k-1}, t_k)(\m_0) \end{pmatrix}.
\end{equation*}\end{linenomath}
If there are $k$ jumps in the time interval $[0,t]$ and if these jumps occur at times $t_1, \dots, t_k$, then $(A^*_t, \M^*_t) = (t-t_k, \theta^k_t(t_1, \dots, t_k)(\m_0))$. For $k=0$, we simply have
\begin{linenomath}\begin{equation*}
    \phi^0_t\left(\begin{pmatrix} a_0\\ \m_0 \end{pmatrix}\right) =
    \begin{pmatrix} a_0+t \\ \theta^0_t(\m_0) \end{pmatrix}.
\end{equation*}\end{linenomath}

If $f$ is bounded (Assumption~\ref{assumption:functions}), we can write $\eta^k(t_1,\dots,t_k;a_0,\m_0)dt_1\dots dt_k$ the probability density over the $k$-first jump times of the process $(A^*_t,M^*_t)_{t\in\R_+}$ having initial condition $(a_0, \m_0)$. We further define the sub-probability density $\nu^k_t(t_1,\dots,t_k;a_0,\m_0)dt_1 \dots dt_k$:
\begin{linenomath}\begin{equation*}
    \nu^k_t(t_1,\dots,t_k;a_0,\m_0) := \mathbbm{1}_{t_k\leq t}\int_t^\infty\eta^{k+1}(t_1,\dots,t_k,t_{k+1};a_0,\m_0)dt_{k+1}.
\end{equation*}\end{linenomath}
Note that the mass of $\nu^k_t$ is the probability of having exactly $k$ jumps in the time interval $[0,t]$. Hence, $\nu^k_t / \int \nu^k_t$ can be interpreted as the probability density over the $k$ jump times knowing that there are exactly $k$ jumps in the time interval $[0,t]$.

Lastly, for all $k \geq 1$, we denote by $\Pi^k_{t_1,\dots,t_k,\m}$ and $\Pi^k_{a,\m}$ the projections
\begin{linenomath}\begin{align*}
    \Pi^k_{t_1,\dots,t_k,\m_0}&: (t_1, \dots, t_k, a_0, \m_0) \mapsto (t_1, \dots, t_k, \m_0), \\
    \Pi^k_{a,\m}&: (t_1, \dots, t_{k-1}, a, \m) \mapsto (a, \m).
\end{align*}\end{linenomath}
By convention, for $k=0$, these projections are the identity. 

We have the path integral representation:

\begin{lemma} \label{lemma:time-marginals}
Grant Assumption~\ref{assumption:functions}. Let $(A^1_t, \M^1_t, x_t)_{t\in\R_+}$ denote the limit process \eqref{eq:limit} for the initial condition $(A^1_0,\M^1_0)\sim u_0$ and $(\E[H_t])_{t\in\R_+}\in \mathcal{C}(\R_+)$. Then, for all $t>0$, the time-marginal $\rho_t:= \mathcal{L}(A^1_t, \M^1_t)$ is given by the representation formula

\begin{linenomath}\begin{equation}\label{eq:path-integral-representation-simple}
    \rho_t = \sum_{k=0}^\infty (\Pi^k_{a,\m}\circ \phi^k_t \circ \Pi^k_{t_1,\dots,t_k,\m_0})_*(\nu^k_t u_0).
\end{equation}\end{linenomath}

If we further grant Assumption~\ref{assumption:gamma_dif} and if $u_0$ is absolutely continuous, then $\rho_t$ is also absolutely continuous.
\end{lemma}

\begin{proof}
Let $\tau_k$ be the time of the $k$-th jump of $\left(A^1_t, \M^1_t, x_t\right)_{t\in \R_+}$. Since $\left(A^1_t, \M^1_t\right)$ is a function of the initial conditions $\left(A^1_0, \M^1_0\right)$ and the jump times $\{\tau_k\}_{k\in\mathbb{N}^*}$, for any continuous and bounded test function $F$ on $\R_+ \times \Rd$, we can write $\E[F(A_t, \M_t)]$ as a `path integral':
\begin{linenomath}\begin{align*}
    &\E[F(A^1_t, \M^1_t)] = \E\left[F(A^1_t, \M^1_t)\mathbbm{1}_{\{t<\tau_1\}}\right] + \sum_{k=1}^\infty \E\left[F(A^1_t, \M^1_t)\mathbbm{1}_{\{\tau_k \leq t < \tau_{k+1}\}}\right] \\
    &= \E\left[F(\phi^0_t(A^1_0,\M^1_0))\mathbbm{1}_{\{t<\tau_1\}}\right] + \sum_{k=1}^\infty \E\left[F(\Pi^k_{a,\m}\circ\phi^k_t(\tau_1, \dots, \tau_k,\M^1_0))\mathbbm{1}_{\{\tau_k \leq t < \tau_{k+1}\}}\right] \\
    &= \int_{\Rd}\int_{\R_+} F(\phi^0_t(a_0,\m_0)) \nu^0_t(a_0,\m_0)u_0(da_0, d\m_0) \\
    &+ \sum_{k=1}^\infty \int_{\Rd}\underbrace{\int_0^t \dots \int_0^t}_{k \text{ times}} F(\Pi^k_{a,\m}\circ\phi^k_t(t_1,\dots,t_k,\m_0))\underbrace{\int_{\R_+}\nu^k_t(t_1,\dots,t_k;a_0,\m_0)u_0(da_0,d\m_0)}_{(\Pi^k_{t_1,\dots,t_k,\m})_*(\nu^k_t u_0)}dt_1\dots dt_k\\
    &= \int_{\Rd}\int_{\R_+} F(a, \m) \biggr(\sum_{k=0}^\infty (\Pi^k_{a,\m}\circ \phi^k_t \circ \Pi^k_{t_1,\dots,t_k,\m_0})_*(\nu^k_t u_0)\biggr)(da,d\m),
\end{align*}\end{linenomath}
whence the representation formula~\eqref{eq:path-integral-representation-simple}.

If $u_0$ is absolutely continuous, then $v^k_tu_0$ is absolutely continuous for all $k\geq0$. If, in addition, Assumption~\ref{assumption:gamma_dif} is granted, then $\phi_t^k$ is a proper local diffeomorphism and $(\Pi^k_{a,\m}\circ \phi^k_t \circ \Pi^k_{t_1,\dots,t_k,\m})_*(\nu^k_t u_0)$ is absolutely continuous for all $k\geq 0$. The probability measure $\rho_t$ is therefore absolutely continuous.
\end{proof}

\subsection{From the limit process to weak solutions}
\begin{proposition} \label{proposition:SDE_to_PDE}
Grant Assumptions~\ref{assumption:functions} and \ref{assumption:gamma_dif}. Further assume that the law of the initial condition $(A^1_0,\M^1_0)$ is the absolutely continuous probability measure $u_0(a,\m)dad\m$ and $(\E[H_t])_{t\in\R_+} = (\bar{H}^1_t)_{t\in\R_+}\in \mathcal{C}(\R_+)$. Then, the time-marginals $(\rho_t)_{t\in\R_+}$ of the law of the process (given by the path integral representation~\eqref{eq:path-integral-representation-simple}) and $(x_t)_{t\in\R_+}$ is a weak solution to~\eqref{eq:PDE}.
\end{proposition}
\begin{proof}
First, we use the path integral representation~\eqref{eq:path-integral-representation-simple} to prove that $\rho \in \mathcal{C}(\R_+, L^1(\R_+ \times \Rd))$. We have to show that for all $T>0$, $\rho \in \mathcal{C}([0,T],L^1(\R_+\times\Rd))$. Let us take $Z^1_t$ from~\eqref{eq:limit_Z}, which counts the number of events in the time interval $[0,t]$. For any $t\in[0,T]$ and any $l\in\mathbb{N}$,
\begin{linenomath}\begin{equation*}
    \sum_{k=l}^\infty\norm{\nu^k_{t} u_0}_{L^1} =\mathbb{P}(Z^1_t\geq l)\leq \mathbb{P}(Z^1_T\geq l).
\end{equation*}\end{linenomath}
Hence, for all $t', t \in [0,T]$,
\begin{linenomath}\begin{multline*}
    \norm{\rho_{t'} - \rho_t}_{L^1} \leq \sum_{k=0}^\infty\norm{(\Pi^k_{a,\m}\circ\phi^k_{t'}\circ\Pi^k_{t_1,\dots,t_k,\m_0})_*(\nu^k_{t'} u_0) - (\Pi^k_{a,\m}\circ\phi^k_t\circ\Pi^k_{t_1,\dots,t_k,\m_0})_*(\nu^k_t u_0)}_{L^1}\\
    \leq \sum_{k=0}^l\norm{(\Pi^k_{a,\m}\circ\phi^k_{t'}\circ\Pi^k_{t_1,\dots,t_k,\m_0})_*(\nu^k_{t'} u_0) - (\Pi^k_{a,\m}\circ\phi^k_t\circ\Pi^k_{t_1,\dots,t_k,\m_0})_*(\nu^k_t u_0)}_{L^1} + 2 \mathbb{P}(Z_T > l).
\end{multline*}\end{linenomath}
Since $\mathbb{P}(Z_T>l)\to 0$ as $l \to \infty$, to show that $\rho \in \mathcal{C}([0,T],L^1(\R_+\times\Rd))$, it suffices to show that for all $k\in\mathbb{N}$, 
\begin{linenomath}\begin{equation*}
   \left((\Pi^k_{a,\m}\circ\phi^k_t\circ\Pi^k_{t_1,\dots,t_k,\m_0})_*(\nu^k_t u_0)\right)_{t\in[0,T]} \in \mathcal{C}([0,T],L^1(\R_+\times\Rd)).
\end{equation*}\end{linenomath}
By the density of $\mathcal{C}_c(\R_+\times\Rd)$ in $L^1(\R_+\times\Rd)$, for any $\epsilon>0$, there exists $\widetilde{u_0}\in\mathcal{C}_c(\R_+\times\Rd)$ such that $\norm{\widetilde{u_0} - u_0}_{L^1}<\frac{\epsilon}{3\norm{f}_\infty^k}$. For all $t\in[0,T]$,
\begin{linenomath}\begin{align*}
    \norm{(\Pi^k_{a,\m}\circ\phi^k_{t}\circ\Pi^k_{t_1,\dots,t_k,\m_0})_*(\nu^k_{t}\widetilde{u_0})- (\Pi^k_{a,\m}\circ\phi^k_{t}\circ\Pi^k_{t_1,\dots,t_k,\m_0})_*(\nu^k_{t}u_0)}_{L_1} &= \norm{\nu^k_t(\widetilde{u_0}-u_0)}_{L^1} \\
    &\leq \norm{f}_\infty^k\norm{\widetilde{u_0} - u_0}_{L^1} \leq \frac{\epsilon}{3}.
\end{align*}\end{linenomath}
Hence, by triangular inequality, it only remains to show that for all $\widetilde{u_0}\in\mathcal{C}_c(\R_+\times\Rd)$,
\begin{linenomath}\begin{equation}
   \left((\Pi^k_{a,\m}\circ\phi^k_t\circ\Pi^k_{t_1,\dots,t_k,\m_0})_*(\nu^k_t \widetilde{u_0})\right)_{t\in[0,T]} \in \mathcal{C}([0,T],L^1(\R_+\times\Rd)). \tag{$\ast$}
\end{equation}\end{linenomath}
Since $\widetilde{u_0}$ is compactly supported, there exists $C>0$ such that $\text{Supp}(\widetilde{u_0})\subset [0,C]\times [-C,+C]^d$. For all $t\in[0,T]$,
\begin{linenomath}\begin{equation*}
    (\Pi^k_{a,\m}\circ\phi^k_t\circ\Pi^k_{t_1,\dots,t_k,\m_0})_*(\nu^k_t \widetilde{u_0}) \leq \mathbbm{1}_{[0,C+T]\times [-C-k\norm{\Gamma}_\infty,+C+k\norm{\Gamma}_\infty]^d}\norm{\widetilde{u_0}}_\infty \in L^1(\R_+\times\Rd).
\end{equation*}\end{linenomath}
Therefore, $(\ast)$ is verified by dominated convergence. This achieves the proof that $\rho \in \mathcal{C}(\R_+, L^1(\R_+ \times \Rd))$.

Now, we verify that $\rho$ satisfies~\eqref{eq:weak} for all test functions. For any $G \in \mathcal{C}^\infty_c(\R_+\times\R_+\times\Rd)$ and any $T>0$, by Itô's formula for jump processes,
\begin{linenomath}\begin{multline*}
    G(T,A^1_T, \M^1_T) = G(0,A^1_0, \M^1_0) + \int_0^T[\partial_t + \partial_a + b(\M^1_t)\cdot \nabla]G(t,A^1_t,\M^1_t)dt \\
    + \int_{[0,T]\times\R_+}\left(G(t,0,\gamma(\M^1_{t-})) - G(t,A^1_{t-}, \M^1_{t-})\right)\mathbbm{1}_{z\leq f(A^1_{t-},\M^1_{t-},x_t)}\pi^1(dt,dz).
\end{multline*}\end{linenomath}
Taking the expectation,
\begin{linenomath}\begin{multline*}
    \E[G(T,A^1_T, \M^1_T)] = \E[G(0,A^1_0, \M^1_0)] + \int_0^T\E\left[[\partial_t + \partial_a + b(\M^1_t)\cdot \nabla]G(t,A^1_t,\M^1_t)\right]dt \\
    + \int_0^T\E\left[\left(G(t,0,\gamma(\M^1_t)) - G(t,A^1_t, \M^1_t)\right)f(A^1_t,\M^1_t,x_t)\right]dt,
\end{multline*}\end{linenomath}
which is equivalent to
\begin{linenomath}\begin{multline}
    \int_{\R_+}\int_{\Rd}G(T,a, \m)\rho_t(a,\m)dad\m \\
    = \int_{\R_+}\int_{\Rd}G(0,a, \m)u_0(a,\m)dad\m + \int_0^T\int_{\R_+}\int_{\Rd}\Big\{[\partial_t + \partial_a + b(\m)\cdot \nabla]G(t,a,\m) \\
    + \left(G(t,0,\gamma(\m)) - G(t,a, \m)\right)f(a,\m,x_t)\Big\}\rho_t(a,\m)dad\m dt. \label{eq:weak_T}
\end{multline}\end{linenomath}
Since $G$ in compactly supported, the $T\to\infty$ limit of~\eqref{eq:weak_T} is \eqref{eq:weak}. This concludes the proof.
\end{proof}

\subsection{Uniqueness of weak solutions}
\begin{proposition}
Grant Assumptions~\ref{assumption:functions} and \ref{assumption:gamma_dif}. For any $(u_0, \bar{H})\in L^1(\R_+\times\R^d, \R_+)\times\mathcal{C}(\R_+)$, the solution to \eqref{eq:PDE} is unique. 
\end{proposition}
\begin{proof}
Let $(\rho, x)$ be a weak solution for some $(u_0, \bar{H})\in L^1(\R_+\times\R^d, \R_+)\times\mathcal{C}(\R_+)$. By Assumption~\ref{assumption:gamma_dif}, the border condition~\eqref{eq:PDE_border} can be written \eqref{eq:border_explicit} and the function
\begin{linenomath}
    \begin{equation}\label{eq:PDE_border_unique}
        (t,\m)\mapsto\mathbbm{1}_{\gamma(\Rd)}(\m)\sum_{\m'\in\gamma^{-1}(\m)}\frac{1}{|\det(\mathbf{J}_\gamma(\m'))|}\int_{\R_+}f(a,\m',x_t)\rho_t(a,\m')da=:p_t(\m)
    \end{equation}
\end{linenomath}
is in $\mathcal{C}(\R_+, L^1(\Rd))$ since $f$ is bounded and Lipschitz with respect to the third variable.

By the standard theory of transport equations with initial datum in $L^1$ (see \cite{Per07}) and treating~\eqref{eq:PDE_border_unique} as a source term, $\rho$ solves
\begin{linenomath}\begin{align}
    &\rho_t(a,\mathbf{m}) \nonumber\\
    &\quad= \begin{cases} 
    u_0(a-t,B_t^{-1}(\m))\exp\left(\int_0^t (\nabla \cdot b)(B_{t-s}^{-1}(\m)) - f(a-t+s,B_{t-s}^{-1}(\m),x_s)ds\right) &\text{if } a \geq t, \\
    p_{t-a}(B_a^{-1}(\m))\exp\left(\int_{t-a}^t (\nabla \cdot b)(B_{t-s}^{-1}(\m)) - f(a-t+s,B_{t-s}^{-1}(\m),x_s)ds\right) &\text{if } 0 < a < t.
    \end{cases}\label{eq:foliation_lagrangian}
\end{align}\end{linenomath}

Using \eqref{eq:PDE_border_unique} and \eqref{eq:foliation_lagrangian}, we have the rough bound $\norm{\rho_t}_{L^1} \leq \norm{u_0}_{L^1}\exp(t \norm{f}_\infty)$:
\begin{linenomath}\begin{align*}
    \norm{\rho_t}_{L^1}&\leq \norm{u_0}_{L^1} + \int_0^t \int_{\Rd}p_{t-a}(\m)d\m da = \norm{u_0}_{L^1} + \int_0^t \int_{\Rd}p_{s}(\m)d\m ds \\
    &\leq \norm{u_0}_{L^1} + \norm{f}_\infty\int_0^t \int_{\Rd}\int_{\R_+}\rho_s(a,\m)da d\m ds
\end{align*}\end{linenomath}
and the bound is obtained using Grönwall's lemma.

Let $(\rho^*, x^*)$ be another weak solution to \eqref{eq:PDE} for the same $(u_0, \bar{H})$. In the following, we derive bounds on the distance $\norm{\rho_t - \rho^*_t}_{L^1} + |x_t - x_t^*|$ and apply Grönwall's lemma. This is relatively straightforward since $f$ and $h$ are bounded and $f$ is Lipschitz with respect to the third variable. For all finite time $T>0$ and for all $t\in[0,T]$,
\begin{linenomath}\begin{align*}
    &\norm{\rho_t - \rho^*_t}_{L^1} \\
    &\quad\leq \int_{\Rd}\int_t^\infty \biggr|u_0(a-t,B^{-1}_t(\m))\exp\left(\int_0^t (\nabla \cdot b)(B_{t-s}^{-1}(\m)) - f(a-t+s,B^{-1}_{t-s}(\m),x_s)ds\right) \\
    &\quad\qquad\quad - u_0(a-t,B^{-1}_t(\m))\exp\left(\int_0^t (\nabla \cdot b)(B_{t-s}^{-1}(\m)) - f(a-t+s,B^{-1}_{t-s}(\m),x^*_s)ds\right)\biggr|dad\m \\
    &\quad\quad+\int_{\Rd}\int_0^t\biggr|  p_{t-a}(B^{-1}_a(\m))\exp\left(\int_{t-a}^t (\nabla \cdot b)(B_{t-s}^{-1}(\m)) - f(a-t+s,B^{-1}_{t-s}(\m),x_s)ds\right) \\
    &\quad\qquad\quad - p^*_{t-a}(B^{-1}_a(\m))\exp\left(\int_{t-a}^t (\nabla \cdot b)(B_{t-s}^{-1}(\m)) - f(a-t+s,B^{-1}_{t-s}(\m),x^*_s)ds\right)\biggr|dad\m \\
    &\quad=: Q_1 + Q_2.
\end{align*}\end{linenomath}
But $Q_1 \leq \norm{u_0}L_f \int_0^t |x_s - x^*_s|ds$, and by triangular inequality (using the shorthand $f(x_t):=f(a,\m,x_t))$),
\begin{linenomath}\begin{align*}
    Q_2 &\leq \left(\int_0^t\int_{\Rd}\rho_{s}(0,\m)d\m ds\right) L_f\int_0^t |x_s - x^*_s|ds + \int_0^t \int_{\Rd} \left|p_{s}(\m) - p^*_{s}(\m)\right|d\m ds \\
    &= \left(\int_0^t\int_{\Rd}\int_{\R_+}f(x_t)\rho_{s}\,dad\m ds\right) L_f\int_0^t |x_s - x^*_s|ds + \int_0^t \int_{\Rd} \int_{\R_+}\left|f(x_t)\rho_{s} - f(x^*_t)\rho^*_{s}\right|\,da d\m ds \\
    &\leq \norm{f}_\infty\int_0^t\norm{\rho_s}_{L^1}ds\;L_f\int_0^t |x_s - x^*_s|ds + \sup_{s\in[0,t]}\norm{\rho_s}_{L^1}L_f\int_0^t|x_s - x^*_s|ds + \norm{f}_\infty\int_0^t\norm{\rho_s - \rho^*_s}ds.
\end{align*}\end{linenomath}
By the rough bound on $\norm{\rho_t}_{L^1}$ established above, we get
\begin{linenomath}\begin{equation*}
    \norm{\rho_t - \rho^*_t}_{L^1} \leq C_{T,0}\int_0^t\norm{\rho_s - \rho^*_s}_{L_1} + |x_s - x^*_s|ds.
\end{equation*}
\end{linenomath}
On the other hand,
\begin{linenomath}\begin{equation*}
    |x_t - x^*_t| \leq \norm{h}_\infty \int_0^t \int_{\Rd} \int_{\R_+}\left|f(x_s)\rho_{s} - f(x^*_s)\rho^*_{s}\right|\,da d\m ds,
\end{equation*}\end{linenomath}
which can be bounded as shown above. Whence,
\begin{linenomath}\begin{equation*}
    \norm{\rho_t - \rho^*_t}_{L^1} + |x_t - x^*_t| \leq C_{T,0}\int_0^t\norm{\rho_s - \rho^*_s}_{L_1} + |x_s - x^*_s|ds.
\end{equation*}
\end{linenomath}
By Grönwall's lemma, $\norm{\rho_t - \rho^*_t}_{L^1} + |x_t - x^*_t| = 0$ for all $t \in [0, T]$. Since this is true for all $T>0$, $(\rho,x) = (\rho^*,x^*)$, which concludes the proof.
\end{proof}

\section*{Acknowledgements}
 I would like to warmly thank Eva L\"{o}cherbach and Wulfram Gerstner for supervising this work and for their comments on this manuscript. I would also like to thank Pablo Ferrari and Monia Capanna for discussions at the initial stages of this project, Claudia Fonte and Stéphane Mischler for discussions on the PDE aspects of this work and Victor Panaretos for general comments on the manuscript. Finally, I would like to thank two anonymous referees whose numerous comments have helped to significantly improve both the form and the content of this manuscript. This research was supported by the Swiss National Science Foundation (grant no.~$200020\_184615$).

\bibliographystyle{plain} 
\bibliography{mybib}

\begin{thebibliography}{10}

\bibitem{BacMas15}
Emmanuel Bacry, Iacopo Mastromatteo, and Jean-Fran{\c{c}}ois Muzy.
\newblock Hawkes processes in finance.
\newblock {\em Market Microstructure and Liquidity}, 1(01):1550005, 2015.

\bibitem{BreMas96}
Pierre Br\'{e}maud and Laurent Massouli\'{e}.
\newblock Stability of nonlinear {H}awkes processes.
\newblock {\em Ann. Probab.}, 24(3):1563--1588, 1996.

\bibitem{CanCar13}
Jos\'{e}~A. Ca\~{n}izo, Jos\'{e}~A. Carrillo, and S\'{\i}lvia Cuadrado.
\newblock Measure solutions for some models in population dynamics.
\newblock {\em Acta Appl. Math.}, 123:141--156, 2013.

\bibitem{Che17}
Julien Chevallier.
\newblock Mean-field limit of generalized {H}awkes processes.
\newblock {\em Stochastic Process. Appl.}, 127(12):3870--3912, 2017.

\bibitem{CorDes13}
Jesus~M Cortes, Mathieu Desroches, Serafim Rodrigues, Romain Veltz, Miguel~A
  Mu{\~n}oz, and Terrence~J Sejnowski.
\newblock Short-term synaptic plasticity in the deterministic
  {T}sodyks--{M}arkram model leads to unpredictable network dynamics.
\newblock {\em Proc. Natl. Acad. Sci. USA}, 110(41):16610--16615, 2013.

\bibitem{CraSor08}
Riley Crane and Didier Sornette.
\newblock Robust dynamic classes revealed by measuring the response function of
  a social system.
\newblock {\em Proc. Natl. Acad. Sci. USA}, 105(41):15649--15653, 2008.

\bibitem{DekLai19}
Marc De~Kamps, Mikkel Lepper{\o}d, and Yi~Ming Lai.
\newblock Computational geometry for modeling neural populations: From
  visualization to simulation.
\newblock {\em PLoS Comput. Biol.}, 15(3):e1006729, 2019.

\bibitem{DemGal15}
Anna De~Masi, Antonio Galves, Eva L\"{o}cherbach, and Errico Presutti.
\newblock Hydrodynamic limit for interacting neurons.
\newblock {\em J. Stat. Phys.}, 158(4):866--902, 2015.

\bibitem{DelFou16}
Sylvain Delattre, Nicolas Fournier, and Marc Hoffmann.
\newblock Hawkes processes on large networks.
\newblock {\em Ann. Appl. Probab.}, 26(1):216--261, 2016.

\bibitem{DitLoe17}
Susanne Ditlevsen and Eva L\"{o}cherbach.
\newblock Multi-class oscillating systems of interacting neurons.
\newblock {\em Stochastic Process. Appl.}, 127(6):1840--1869, 2017.

\bibitem{DuaLoe19}
Aline Duarte, Eva L\"{o}cherbach, and Guilherme Ost.
\newblock Stability, convergence to equilibrium and simulation of non-linear
  {H}awkes processes with memory kernels given by the sum of {E}rlang kernels.
\newblock {\em ESAIM Probab. Stat.}, 23:770--796, 2019.

\bibitem{FonSch21}
Claudia Fonte and Valentin Schmutz.
\newblock Long time behavior of an age and leaky memory-structured neuronal
  population equation.
\newblock {\em arXiv preprint arXiv:2106.11110}, 2021.

\bibitem{FouGui15}
Nicolas Fournier and Arnaud Guillin.
\newblock On the rate of convergence in {W}asserstein distance of the empirical
  measure.
\newblock {\em Probab. Theory Related Fields}, 162(3-4):707--738, 2015.

\bibitem{FouLoe16}
Nicolas Fournier and Eva L{\"o}cherbach.
\newblock On a toy model of interacting neurons.
\newblock In {\em Annales de l'Institut Henri Poincar{\'e}, Probabilit{\'e}s et
  Statistiques}, volume~52, pages 1844--1876. Institut Henri Poincar{\'e},
  2016.

\bibitem{GalLoe16}
Antonio Galves and Eva L\"{o}cherbach.
\newblock Modeling networks of spiking neurons as interacting processes with
  memory of variable length.
\newblock {\em J. SFdS}, 157(1):17--32, 2016.

\bibitem{GalLoe20}
Antonio Galves, Eva L\"{o}cherbach, Christophe Pouzat, and Errico Presutti.
\newblock A system of interacting neurons with short term synaptic
  facilitation.
\newblock {\em J. Stat. Phys.}, 178(4):869--892, 2020.

\bibitem{GerEde13}
Felipe Gerhard, Tilman Kispersky, Gabrielle~J Gutierrez, Eve Marder, Mark
  Kramer, and Uri Eden.
\newblock Successful reconstruction of a physiological circuit with known
  connectivity from spiking activity alone.
\newblock {\em PLoS Comput. Biol.}, 9(7):e1003138, 2013.

\bibitem{Ger95}
Wulfram Gerstner.
\newblock Time structure of the activity in neural network models.
\newblock {\em Phys. Rev. E}, 51(1):738, 1995.

\bibitem{Ger00}
Wulfram Gerstner.
\newblock Population dynamics of spiking neurons: fast transients, asynchronous
  states, and locking.
\newblock {\em Neural Comput.}, 12(1):43--89, 2000.

\bibitem{GerKis02}
Wulfram Gerstner and Werner~M Kistler.
\newblock {\em Spiking neuron models: Single neurons, populations, plasticity}.
\newblock Cambridge university press, 2002.

\bibitem{GerKis14}
Wulfram Gerstner, Werner~M Kistler, Richard Naud, and Liam Paninski.
\newblock {\em Neuronal dynamics: From single neurons to networks and models of
  cognition}.
\newblock Cambridge University Press, 2014.

\bibitem{GerHem92}
Wulfram Gerstner and J~Leo van Hemmen.
\newblock Associative memory in a network of `spiking'neurons.
\newblock {\em Netw. Comput. Neural Syst.}, 3(2):139--164, 1992.

\bibitem{Haw71}
Alan~G. Hawkes.
\newblock Spectra of some self-exciting and mutually exciting point processes.
\newblock {\em Biometrika}, 58:83--90, 1971.

\bibitem{Haw18}
Alan~G. Hawkes.
\newblock Hawkes processes and their applications to finance: a review.
\newblock {\em Quant. Finance}, 18(2):193--198, 2018.

\bibitem{JolRau06}
Renaud Jolivet, Alexander Rauch, Hans-Rudolf L{\"u}scher, and Wulfram Gerstner.
\newblock Predicting spike timing of neocortical pyramidal neurons by simple
  threshold models.
\newblock {\em Journal of computational neuroscience}, 21(1):35--49, 2006.

\bibitem{KipLan99}
Claude Kipnis and Claudio Landim.
\newblock {\em Scaling limits of interacting particle systems}, volume 320 of
  {\em Fundamental Principles of Mathematical Sciences}.
\newblock Springer-Verlag, Berlin, 1999.

\bibitem{LaiDek17}
Yi~Ming Lai and Marc de~Kamps.
\newblock Population density equations for stochastic processes with memory
  kernels.
\newblock {\em Phys. Rev. E}, 95(6):062125, 2017.

\bibitem{LamTul18}
Regis~C Lambert, Christine Tuleau-Malot, Thomas Bessaih, Vincent Rivoirard,
  Yann Bouret, Nathalie Leresche, and Patricia Reynaud-Bouret.
\newblock Reconstructing the functional connectivity of multiple spike trains
  using {H}awkes models.
\newblock {\em J. Neurosci. Methods}, 297:9--21, 2018.

\bibitem{Loe17}
Eva L\"{o}cherbach.
\newblock Spiking neurons: interacting {H}awkes processes, mean field limits
  and oscillations.
\newblock In {\em Journ\'{e}es {MAS} 2016 de la {SMAI}---{P}h\'{e}nom\`enes
  complexes et h\'{e}t\'{e}rog\`enes}, volume~60 of {\em ESAIM Proc. Surveys},
  pages 90--103. EDP Sci., Les Ulis, 2017.

\bibitem{MehSch20}
Sima Mehri, Michael Scheutzow, Wilhelm Stannat, and Bian~Z. Zangeneh.
\newblock Propagation of chaos for stochastic spatially structured neuronal
  networks with delay driven by jump diffusions.
\newblock {\em Ann. Appl. Probab.}, 30(1):175--207, 2020.

\bibitem{Mel96}
Sylvie M\'{e}l\'{e}ard.
\newblock Asymptotic behaviour of some interacting particle systems;
  {M}c{K}ean-{V}lasov and {B}oltzmann models.
\newblock In {\em Probabilistic models for nonlinear partial differential
  equations ({M}ontecatini {T}erme, 1995)}, volume 1627 of {\em Lecture Notes
  in Math.}, pages 42--95. Springer, Berlin, 1996.

\bibitem{MonBar08}
Gianluigi Mongillo, Omri Barak, and Misha Tsodyks.
\newblock Synaptic theory of working memory.
\newblock {\em Science}, 319(5869):1543--1546, 2008.

\bibitem{MulBue07}
Eilif Muller, Lars Buesing, Johannes Schemmel, and Karlheinz Meier.
\newblock Spike-frequency adapting neural ensembles: beyond mean adaptation and
  renewal theories.
\newblock {\em Neural Comput.}, 19(11):2958--3010, 2007.

\bibitem{NauGer12}
Richard Naud and Wulfram Gerstner.
\newblock Coding and decoding with adapting neurons: a population approach to
  the peri-stimulus time histogram.
\newblock {\em PLoS Comput. Biol.}, 8(10):e1002711, 2012.

\bibitem{Oga99}
Yosihiko Ogata.
\newblock Seismicity analysis through point-process modeling: A review.
\newblock In {\em Seismicity patterns, their statistical significance and
  physical meaning}, pages 471--507. Springer, 1999.

\bibitem{PakPer10}
Khashayar Pakdaman, Beno\^{\i}t Perthame, and Delphine Salort.
\newblock Dynamics of a structured neuron population.
\newblock {\em Nonlinearity}, 23(1):55--75, 2010.

\bibitem{PakPer13}
Khashayar Pakdaman, Beno\^{\i}t Perthame, and Delphine Salort.
\newblock Relaxation and self-sustained oscillations in the time elapsed neuron
  network model.
\newblock {\em SIAM J. Appl. Math.}, 73(3):1260--1279, 2013.

\bibitem{PakPer14}
Khashayar Pakdaman, Beno\^{\i}t Perthame, and Delphine Salort.
\newblock Adaptation and fatigue model for neuron networks and large time
  asymptotics in a nonlinear fragmentation equation.
\newblock {\em J. Math. Neurosci.}, 4:Art. 14, 26, 2014.

\bibitem{PayGue21}
Alexandre Payeur, Jordan Guerguiev, Friedemann Zenke, Blake~A Richards, and
  Richard Naud.
\newblock Burst-dependent synaptic plasticity can coordinate learning in
  hierarchical circuits.
\newblock {\em Nat. Neurosci.}, 24(7):1010--1019, 2021.

\bibitem{Per07}
Beno\^{\i}t Perthame.
\newblock {\em Transport equations in biology}.
\newblock Frontiers in Mathematics. Birkh\"{a}user Verlag, Basel, 2007.

\bibitem{PilShl08}
Jonathan~W Pillow, Jonathon Shlens, Liam Paninski, Alexander Sher, Alan~M
  Litke, EJ~Chichilnisky, and Eero~P Simoncelli.
\newblock Spatio-temporal correlations and visual signalling in a complete
  neuronal population.
\newblock {\em Nature}, 454(7207):995--999, 2008.

\bibitem{Qui16}
Crist\'{o}bal Qui\~{n}inao.
\newblock A microscopic spiking neuronal network for the age-structured model.
\newblock {\em Acta Appl. Math.}, 146:29--55, 2016.

\bibitem{RaaDit20}
Mads~Bonde Raad, Susanne Ditlevsen, and Eva L\"{o}cherbach.
\newblock Stability and mean-field limits of age dependent {H}awkes processes.
\newblock {\em Ann. Inst. Henri Poincar\'{e} Probab. Stat.}, 56(3):1958--1990,
  2020.

\bibitem{ReyRiv13}
Patricia Reynaud-Bouret, Vincent Rivoirard, and Christine Tuleau-Malot.
\newblock Inference of functional connectivity in neurosciences via {H}awkes
  processes.
\newblock In {\em 2013 IEEE Glob. Conf. Signal Inf. Process.}, pages 317--320.
  IEEE, 2013.

\bibitem{ReySch10}
Patricia Reynaud-Bouret and Sophie Schbath.
\newblock Adaptive estimation for {H}awkes processes; application to genome
  analysis.
\newblock {\em Ann. Statist.}, 38(5):2781--2822, 2010.

\bibitem{SchChi19}
Tilo Schwalger and Anton~V Chizhov.
\newblock Mind the last spike---firing rate models for mesoscopic populations
  of spiking neurons.
\newblock {\em Curr. Opin. Neurobiol.}, 58:155--166, 2019.

\bibitem{Szn91}
Alain-Sol Sznitman.
\newblock Topics in propagation of chaos.
\newblock In {\em \'{E}cole d'\'{E}t\'{e} de {P}robabilit\'{e}s de
  {S}aint-{F}lour {XIX}---1989}, volume 1464 of {\em Lecture Notes in Math.},
  pages 165--251. Springer, Berlin, 1991.

\bibitem{ToyRad09}
Taro Toyoizumi, Kamiar~Rahnama Rad, and Liam Paninski.
\newblock Mean-field approximations for coupled populations of generalized
  linear model spiking neurons with {M}arkov refractoriness.
\newblock {\em Neural Comput.}, 21(5):1203--1243, 2009.

\bibitem{Tru16}
Wilson Truccolo.
\newblock From point process observations to collective neural dynamics:
  {N}onlinear {H}awkes process {GLM}s, low-dimensional dynamics and coarse
  graining.
\newblock {\em J. Physiol. Paris}, 110(4):336--347, 2016.

\bibitem{TruEde05}
Wilson Truccolo, Uri~T Eden, Matthew~R Fellows, John~P Donoghue, and Emery~N
  Brown.
\newblock A point process framework for relating neural spiking activity to
  spiking history, neural ensemble, and extrinsic covariate effects.
\newblock {\em J. Neurophysiol.}, 93(2):1074--1089, 2005.

\bibitem{TsoPaw98}
Misha Tsodyks, Klaus Pawelzik, and Henry Markram.
\newblock Neural networks with dynamic synapses.
\newblock {\em Neural Comput.}, 10(4):821--835, 1998.

\bibitem{WilCow72}
Hugh~R Wilson and Jack~D Cowan.
\newblock Excitatory and inhibitory interactions in localized populations of
  model neurons.
\newblock {\em Biophys. J.}, 12(1):1--24, 1972.

\bibitem{ZucReg02}
Robert~S Zucker and Wade~G Regehr.
\newblock Short-term synaptic plasticity.
\newblock {\em Annu. Rev. Physiol.}, 64(1):355--405, 2002.

\end{thebibliography}

\end{document}